\documentclass{paper}
\usepackage[T1]{fontenc}
\usepackage[utf8]{inputenc}
\usepackage{makeidx}
\usepackage{graphicx}
\usepackage{float} 
\usepackage{amssymb}
\usepackage{amsmath}
\usepackage{latexsym}
\usepackage{wasysym}
\usepackage{amsthm}
\usepackage{natbib}
\begin{document}

%\markboth{Author's Name}{Paper Title}
\title{Numerical convergence of the block-maxima approach to the Generalized Extreme Value  distribution}

%\author{D. Faranda, V. Lucarini, G. Turchetti, S. Vaienti, j. Wouters}
\author{Faranda, Davide\\
{\small Department of Mathematics and Statistics, University of Reading;}\\
{\small Whiteknights, PO Box 220, Reading RG6 6AX, UK. d.faranda@pgr.reading.ac.uk}\\
\\
Lucarini, Valerio\\
{\small Department of Meteorology, University of Reading;}\\
{\small Department of Mathematics and Statistics, University of Reading;}\\
{\small Whiteknights, PO Box 220, Reading RG6 6AX, UK. v.lucarini@reading.ac.uk}\\
\\
Turchetti, Giorgio\\
{\small Department of Physics, University of Bologna.INFN-Bologna}\\
{\small Via Irnerio 46, Bologna, 40126, Italy. turchett@bo.infn.it}\\
\\
Vaienti, Sandro\\
{\small UMR-6207, Centre de Physique Th\'eorique, CNRS, Universit\'es d'Aix-Marseille I,II,}\\
{\small Universit\'e du Sud Toulon-Var and FRUMAM}\\
{\small (F\'ed\'eration de Recherche des Unit\'es de Math\'ematiques de Marseille);}\\
{\small CPT, Luminy, Case 907, 13288 Marseille Cedex 09, France.}\\
{\small vaienti@cpt.univ-mrs.fr}\\
}
\date{}
\maketitle

\begin{abstract}

In this paper we perform an analytical and numerical study of Extreme Value distributions in discrete dynamical systems. In this setting, recent works   have shown how to get a statistics of extremes in agreement with the classical Extreme Value Theory. We pursue these investigations by   giving  analytical expressions  of Extreme Value distribution parameters  for maps that have an absolutely continuous invariant measure. We compare these analytical results with numerical experiments in which we study the convergence to  limiting distributions using the so called block-maxima approach, pointing out in which cases we obtain robust estimation of parameters. In regular maps for which mixing properties do  not  hold, we show that 
the fitting procedure to the classical Extreme Value Distribution fails, as expected. However, we obtain an empirical distribution that can be explained starting from a different observable function for which  \citet{nicolis} have found analytical results.

\end{abstract}

\section{Introduction}
 Extreme Value Theory (EVT) was first developed  by   \citet{fisher} and formalized by \citet{gnedenko}  which showed that the distribution of the block-maxima of a sample of independent identically distributed (i.i.d) variables converges to a member of the so-called Extreme Value (EV) distribution. It arises from the study of stochastical series that is of great interest in different disciplines: it has been applied  to extreme floods \citep{gumbel1941}, \citep{sv}, \citep{friederichs}, amounts of large insurance losses \citep{brodin}, \citep{cruz2002}; extreme earthquakes \citep{sornette}, \citep{cornell}, \citep{burton1979seismic}; meteorological  and climate events \citep{felici1},  \citep{felici2},\citep{vitolo},  \citep{alttresh}, \citep{nicholis1997clivar}, \citep{smith1989extreme}.  All these events have a relevant impact on socioeconomic activities and it is crucial to find a way to understand and, if possible, forecast them \citep{hallerberg2008influence}, \citep{kantz2006dynamical}.\\
The attention of the scientific community to the problem of modeling extreme values is growing. An extensive account of recent results and relevant applications is given in \citet{ghil2010extreme}. Such an interest is mainly due to the fact that this theory is also important in defining risk factor in a wide class of applications  such as the modeling of  financial risk after the  significant instabilities in financial markets worldwide \citep{gilli2006application}, \citep{longin2000value}, \citep{embrechts1999extreme},   the analysis of seismic  and hydrological risk    \citep{burton1979seismic}, \citep{martins2000generalized}. Even if  the probability of extreme events decreases with their magnitude, the damage that they may bring increases rapidly with the magnitude as does the cost of protection against them \citet{nicolis}.\\
From a theoretical point of view, extreme values represent extreme fluctuations of a system. Very recently, many authors have shown clearly how the statistics of global observables in correlated systems can be related to EV statistics \citep{dahlstedt2001universal}, \citep{bertin2005global}. \citet{clusel2008global} have shown how to connect fluctuations of global additive quantities, like total energy or magnetization , by statistics of sums of  random variables in such a way that it is possible to  identify a class of  random variables whose sum follows an extreme value distributions.\\
The so called block-maxima approach is widely used in EVT since it represents a very natural way to look at extremes. It consists of dividing the data series of some observable into bins of equal length and selecting the maximum (or the minimum) value in each of them \citep{coles}.
When dealing with climatological or financial data, since we usually have  limited data-set, the main problem in applying EVT is related to the choice of a sufficiently large statistics of extremes provided that each bin contains a  suitable number of observations. Therefore a smart balance between number of maxima and  observations per bin is needed \citep{felici1}, \citep{katz}, \citep{katz2}, \citep{katz2005statistics}.\\
Recently a number of alternative approaches have been studied. One consists in looking at exceedance over high thresholds rather than maxima over fixed time periods. While the idea of looking at extreme value problems from this point of view is very old, the development of a modern theory has started with \citet{todorovic1970stochastic} that have proposed the so called Peaks Over Threshold approach. At the same time there was a mathematical development of procedures based on a certain number of extreme order statistics \citep{pickands1975statistical}, \citep{hill1975simple} and the Generalized Pareto distribution for excesses over thresholds \citep{smith1984threshold}, \citep{davison1984modelling}, \citep{davison1990models}. \\

Since dynamical systems theory can be used to understand features of physical systems like climate and forecast financial behaviors, many authors have studied how to extend EVT to these field. When dealing with dynamical systems we have to know what kind of properties (i.e. stability, degree of mixing, correlations decay) are related to Gnedenko's hypotheses and also which observables we  must consider in order to obtain an EV distribution. Furthermore, even if the convergence is achieved, we should evaluate how fast it is depending on all parameters and properties used.  Empirical studies show that in some cases a dynamical observable obeys to the extreme value statistics even if the convergence is highly dependent on the kind of observable we choose \citep{vannitsem}, \citep{vitolo}, \citep{vitolo2009robust}. For example,  \citet{balakrishnan} and more recently \citet{nicolis} and \citet{haiman} have shown that for regular orbits of dynamical systems we don't expect to find convergence to EV distribution.\\
The first rigorous  mathematical approach to extreme value theory in dynamical systems goes back to the pioneer paper by P. Collet in 2001 \citep{collet2001statistics}. Collet got the Gumbel Extreme Value Law (see below) for certain one-dimensional non-uniformly
hyperbolic maps which admit an absolutely continuous invariant measure and exhibit exponential decay of correlations. Collet's approach used Young towers \citep{young}, \citep{young2}  and his suggestion was successively applied to other systems. Before quoting them, we would like to point out that Collet was able to establish a few conditions (usually called $D$ and $D'$) and which have been introduced by Leadbetter \citep{leadbetter} with the aim to associate to the stationary stochastic process given by the dynamical system, a new stationary independent sequence which enjoyed one of the classical three extreme value laws, and this law could be pulled back to the original dynamical sequence. Conditions $D$ and $D'$ require a sort of independence of the stochastic dynamical sequence in terms of uniform mixing condition on the distribution functions. Condition $D$ was successively  improved by Freitas and Freitas \citep{freitas2008}, in the sense that they introduced a new condition, called $D_2$, which is weaker than $D$ and that could be checked directly by estimating the rate of decay of correlations for H\"older observables \footnote{We briefly state here the two conditions, we defer to the next section for more details about the quantities introduced. If $X_n, n\ge 0$ is a stochastic process, we define  $M_{j,l}\equiv \{X_j, X_{j+1},\cdots,X_{j+l}\}$ and we put $M_{0,m}=M_m$. Moreover we set  $a_m$ and $b_m$ two normalising sequences and  $u_m=x/a_m +b_m$, where $x$ is a real number, cf. next section for the meaning of these variables. The condition $D_2(u_m)$ holds for the sequence $X_m$ if for any integer $l, t,m$ we have $|\nu (X_0>u_m, M_{t,l}\le u_m)-\nu (X_0>u_m)\nu (M_{t,l}\le u_m)|\le \gamma(m,t)$, where $\gamma(m,t)$ is non-increasing in $t$ for each $m$ and $m\gamma(m,t_m)\rightarrow 0$ as $m\rightarrow \infty$ for some sequence $t_m=o(m)$, $t_m \rightarrow \infty$.\\ We say condition $D'( u_m)$ holds for the sequence $X_m$ if $\lim_{k\rightarrow \infty}\limsup_{m} m\sum_{j=1}^{[m/k]}\nu(X_0>u_m, X_j>u_m)=0$. Whenever the process is given by the iteration of a dynamical systems, the previous two conditions could also be formulated in terms of decay of correlation integrals, see \citet{freitas2008link}, \citet{gupta2010extreme}.}. We notice that conditions $D_2$ and $D'$ allow immediately to get Extreme Value Laws for absolutely continuous invariant measures for uniformly one-dimensional expanding dynamical systems: this is the case for instance of the 1-D maps with constant density studied in Sect. 3 below. Another interesting issue of Collet's paper was the choice of the observables $g$'s whose values along the orbit of the dynamical systems constitute the sequence of events upon which we successively search for the partial maximum. Collet considered a function $g(\mbox{dist}(x, \zeta))$ of the distance with respect to a given point $\zeta$,   with the aim that $g$ achieves a global maximum at almost all points $\zeta$ in the phase space; for example $g(x)=-\log x$. Using a different $g$, Freitas and Freitas \citep{moreira} were able to get the Weibull law for the family of quadratic maps with the Benedicks-Carlesson parameters and for $\zeta$ taken as the critical point or  the critical value, so improving the previous results by Collet who did not keep such values in his set of full measure. \\The latter paper \citep{moreira} strongly relies on condition $D_2$; this condition has  also been invoked to establish the extreme value laws on towers which model dynamical systems with stable foliations (hyperbolic billiards, Lozi maps, H\'enon diffeomorphisms, Lorenz maps and flows). This is the content of the paper by Gupta, Holland and Nicol \citep{gupta2009extreme}. We point out  that the observable $g$ was taken in one of three different classes  $g_1, g_2, g_3$, see Sect. 2 below, each one being again a function of the distance with respect to a given point $\zeta$. The choice of these particular forms for the $g$'s is just to fit with the necessary and sufficient condition on the tail of the distribution function $F(u)$, see next section, in order to exist a non-degenerate limit distribution for the partial maxima \citep{freitas}, \citep{holland}. The paper \citet{gupta2009extreme} also covers the  easier case of uniformly hyperbolic diffeomorphisms, for instance the Arnold Cat map which we studied in Sect. 3.2. \\ Another major step in this field   was achieved by establishing a connection between the extreme value laws and the statistics of first return and hitting times, see the papers by Freitas, Freitas and Todd \citep{freitas}, \citep{freitasNuovo}. They showed in particular that for dynamical systems preserving an absolutely continuous invariant measure or a singular continuous invariant measure $\nu$ , the existence of an exponential hitting time statistics on balls around $\nu$ almost any point $\zeta$ implies the existence of extreme value laws for one of the observables of type $g_i, i=1,2,3$ described above. The converse is also true, namely if we have an extreme value law which applies to the observables of type $g_i, i=1,2,3$ achieving a maximum at $\zeta$, then we have exponential hitting time statistics to balls with center $\zeta$. Recently these results have been generalized to local returns around balls centered at periodic points \citep{freitas2010extremal}. We would like to point out that the equivalence between extreme values laws and the hitting time statistics allowed to prove the former for broad classes of systems for which the statistics of recurrence were known, for instance for expanding maps in higher dimension.\\

In this work we consider a few aspects of the extreme value theory applied to dynamical systems throughout both analytical results and numerical experiments. In particular we analyse the convergence to EV limiting distributions pointing out how robust are parameters estimations. Furthermore, we check the consistency of block-maxima approach highlighting  deviations
from theoretical expected behavior depending on the number of maxima and number of block-observation. To perform our analysis we use low dimensional maps with different properties: mixing maps in which we expect to find convergence to EV distributions and regular maps where the convergence is not ensured.\\
The work is organised as follow: in section 2 we briefly recall methods and results of EVT for independent and identical distributed (i.i.d.) variables and dynamical systems. In section 3 we  explicitly compute theoretical expected distributions parameter in respect to the observable functions of type $g_i, i=1,2,3 $ for map that have constant density measure. Numerical experiments on low dimensional maps are presented.  In section 4 we show that it is possible to derive an asymptotic expression of normalising sequences when the density measure is not constant. As an example we derive the explicit expressions for the Logistic map.  Eventually, in section 5 we repeat the experiment for regular maps showing that extreme values laws do not follow from numerical experiments.

\section{Background on EVT}
\label{background}

\citet{gnedenko}  studied the convergence of maxima of i.i.d. variables $$X_0, X_1, .. X_{m-1}$$ with cumulative distribution (cdf) $F(x)$ of the form: 
$$F(x)=P\{a_m(M_m-b_m) \leq x\}$$
Where $a_m$ and $b_m$ are normalising sequences and $M_m=\max\{ X_0,X_1, ..., X_{m-1}\}$. It may be rewritten as $F(u_m)=P\{M_m \leq u_m\}$ where $u_m=x/a_m +b_m$. Such types of normalising sequences converge to one of the three type of Extreme Value (EV) distribution if  necessary and sufficient conditions on parent distribution of $X_i$ variables are satisfied \citep{leadbetter}.
EV distributions include the following three families:
\begin{itemize}
\item Gumbel  distribution  (type 1):
\begin{equation}
F(x)=
\exp{\{-e^{-x}\}}  \quad x \in \mathbb{R}
\label{gumbel}
\end{equation}
\item Fr\'echet  distribution (type 2):
\begin{equation}
\begin{cases}
F(x)=0  &	x \leq 0 \\
F(x)=\exp{\left\{-x^{1/\xi} \right\}} & x > 0\\
\end{cases}
\label{Fr\'echet}
\end{equation}
\item Weibull  distribution (type 3):
\begin{equation}
\begin{cases}
F(x)=\exp{\left\{- \left(-x \right)^{-1/\xi} \right\}} & \quad x < 0 \\
F(x)=1 & \quad x \geq 0 \\
\end{cases}
\label{weibull}
\end{equation}
\end{itemize}

Let us define the right endpoint $x_F$ of a distribution function $F(x)$ as:
\begin{equation}
x_F=\sup\{ x: F(x)<1\}
\end{equation}

then, it is possible to compute normalising sequences   $a_m$ and $b_m$ using the following corollary of Gnedenko's theorem :\\
\textbf{Corollary (Gnedenko):}  \textit{The normalizing sequences $a_m$ and $b_m$ in the convergence of normalized maxima $P\{a_m(M_m - b_m) \leq x\} \to F(x)$ may be taken (in order of increasing complexity) as:}

\begin{itemize}

\item \textit{Type 1:} $\quad a_m=[G(\gamma_m)]^{-1}, \quad b_m=\gamma_m$;

\item \textit{Type 2:} $\quad a_m=\gamma_m^{-1}, \quad b_m=0 \mbox{ or } b_m=c\cdot m^{-\xi}$;

\item \textit{Type 3:} $\quad a_m=(x_F-\gamma_m)^{-1}, \quad b_m=x_F$;

\end{itemize}
\textit{where}
\begin{equation}
\gamma_m=F^{-1}(1-1/m)=\inf\{x; F(x) \geq 1-1/m\}
\label{gamma}
\end{equation}

\begin{equation}
G(t)=\int_t^{x_F} \frac{1-F(u)}{1-F(t)}du, \quad  t<x_F
\label{gneg}
\end{equation}

and $c \in \mathbb{R}$ is a constant. It is important to remark that the choice of normalising sequences is not unique \citep{leadbetter}. For example for $b_m$ of type 2 distribution it is possible to choose
either $b_m=0$ or $b_m=c \cdot m^{-\xi}$. In particular, we will use the last one since it is a more general choice that ensure the convergence for a much broader class of initial distributions \citep{beirlant}.\\
%\paragraph{Generalized Extreme Value Distribution approach.}
 Instead of Gnedenko's approach it is possible to fit unnormalized data directly  to a single family of generalized distribution called GEV distribution with cdf:

\begin{equation}
F_{G}(x; \mu, \sigma, \xi)=\exp\left\{-\left[1+{\xi}\left(\frac{x-\mu}{\sigma}\right)\right]^{-1/{\xi}}\right\}
\label{cumul}
\end{equation}

which holds for $1+{\xi}(x-\mu)/\sigma>0 $, using $\mu \in \mathbb{R}$ (location parameter) and $\sigma>0$ (scale parameter) as scaling constants in place of $b_m$, and $a_m$ \citep{pickands}. ${\xi} \in \mathbb{R}$ is the shape parameter also called the tail index: when ${\xi} \to 0$, the distribution corresponds to a
Gumbel type. When the index is negative, it corresponds to a Weibull; when the index is positive,
it corresponds to a Fr\'echet.\\

In order to adapt the extreme value theory to dynamical systems, we will consider the stationary stochastic process $X_0,X_1,...$  given by:

\begin{equation}
X_m(x)=g(\mbox{dist}(f^m (x), \zeta)) \qquad \forall m \in \mathbb{N}
\label{sss}
\end{equation}

where 'dist' is a Riemannian metric on $\Omega$, $\zeta$ is a given point and $g$ is an observable function, and whose partial maximum is defined as:

\begin{equation}
{M_m}= \max\{ X_0, ... , X_{m-1} \}
\label{maxi}
\end{equation}

The probability measure will be here an invariant measure $\nu$ for the dynamical system. As we anticipated in the Introduction, we will use   three types of observables $g_i,i=1,2,3$,  suitable to obtain one of the three types of EV distribution  for normalised maxima:

\begin{equation}
g_1(x)= -\log(\mbox{dist}(x,\zeta))
\label{g1}
\end{equation}

\begin{equation}
g_2(x)=\mbox{dist}(x, \zeta)^{-1/\alpha}
\label{g2}
\end{equation}

\begin{equation}
g_3(x)=C - \mbox{dist}(x,\zeta)^{1/\alpha}
\label{g3}
\end{equation}

where $C$ is a constant and $\alpha>0 \in \mathbb{R}$.\\ 
These three type of functions are representative of broader classes which are defined, for instance, throughout equations (1.11) to (1.13) in Freitas et al. [2009]; we now explain the reasons and the meaning of these choices. First of all these functions have in common the following properties: (i) they are defined on the positive semi-axis $[0,\infty]$ with values into $\mathbb{R}\cup \{+\infty\}$; (ii) $0$ is a global maximum, possibly equal to $+\infty$; (iii) they are a strictly decreasing bijection in a neighborhood $V$ of $0$ with image $W$. Then we consider three types of behavior which generalize the previous specific choices:\\
{\em Type 1}: there is a strictly positive function $p:W\rightarrow  \mathbb{R}$ such that $\forall y\in \mathbb{R}$ we have $$ \lim_{s\rightarrow g_1(0)}\frac{g_1^{-1}(s+yp(s))}{g_1^{-1}(s)}=e^{-y}
$$
{\em Type 2}: $g_2(0)=+\infty$ and there exists $\beta>0$ such that $\forall y>0$ we have $$ \lim_{s\rightarrow \infty}\frac{g_2^{-1}(sy)}{g_2^{-1}(s)}=y^{-\beta}
$$
{\em Type 3}: $g_3(0)=D<+\infty$ and there exists $\gamma>0$ such that $\forall y>0$ we have $$ \lim_{s\rightarrow 0}\frac{g_3^{-1}(D-sy)}{g_3^{-1}(D-s)}=y^{\gamma}
$$

The Gnedenko corollary says that the different kinds of extreme value laws are determined by the distribution of $F(u)=\nu (X_0\le u)$ and by the right endpoint of $F$, $x_F$. We will see in the next section that the local invertibility of $g_i, i=1,2,3$ in the neighborhood of $0$ together with the Lebesgue's differentiation theorem (which basically says that whenever the measure $\nu$ is absolutely continuous with respect to Lebesgue with density $\rho$, the the measure of a ball $B_{\delta}(x_0)$ of radius $\delta$ centered around almost any point $x_0$ scales like $\delta \rho(x_0)$), allow us to compute the tail of $F$, in fact we have $$  1-F(u) \sim \rho(\zeta)|B_{g^{-1}(u)}(\zeta)|,$$ where $g$ is any of the three types of functions introduced in (10) to (12) and $|A|$ denotes the diameter of the set $A$. As we said above the tail of $F$ determines the three limit laws for partial maximum of i.i.d. sequences. In  particular    Th. 1.6.2. in \citet{leadbetter} specifies  what kind of conditions the distribution function $F$ must verify to get one specific law: the above type 1,2,3 assumptions  are just the translation in terms of the shape of  $g_i$ of  the conditions on the tail of $F$.

\section{Distribution of Extremes in mixing maps with constant density measure}

Our goal is to use a block-maxima approach and fit our unnormalised  data to a GEV distribution; for that  it will be  necessary to find a linkage among $a_m$, $b_m$, $\mu$ and $\sigma$. At this regard we will use Gnedenko's corollary to compute normalising sequences showing that they correspond to the parameter we obtain fitting directly data to GEV distribution.

We derive the correct expression for mixing maps with constant density measure and the asymptotic behavior for logistic map that is a case of non-constant density measure. 

\subsection{Asymptotic sequences}

In this section we will consider the case of uniformly hyperbolic maps which preserve the Lebesgue measure (the density $\rho=1$) and satisfy the conditions $D_2$ and $D'$, sufficient to get extreme valuers distributions. For the second map, the algebraic automorphisms of the torus better known as the Arnold cat map, the existence of extreme value laws follows from the theory developed in \citet{gupta2009extreme}.
Starting from the definitions provided by Gnedenko we derive as a novel result the exact expression for the normalising sequences $a_m$ and $b_m$.

\paragraph{ Case 1: $\mathbf{g_1}$(x)= -log(dist(x,$\mathbf{\zeta}$)).} 

By  equations \ref{sss} and \ref{maxi} we know that:
\begin{equation}
\begin{split}
1-F(u) &=1-\nu(g(\mbox{dist}(x, \zeta))\leq u) \\
      &= 1 - \nu( -\log( \mbox{dist}(x, \zeta)) \leq u) \\
      &= 1 - \nu( \mbox{dist}(x, \zeta) \geq e^{-u})\\
\end{split}
\label{g1_e1}
\end{equation}

and the last line is justified by using Lebesgue's Differentiation Theorem. Then, for maps with constant density measure, we can write:

\begin{equation}
1- F(u)  \simeq \nu( B_{e^{-u}}(\zeta)) = \Omega_de^{-ud}\\
\label{g1_e2}
\end{equation}

where $d$ is the dimension of the space and $\Omega_d$ is a constant.
To use Gnedenko corollary it is necessary to calculate $u_F$

$$u_F = \sup\{ u ; F(u) < 1\}$$

in this case $u_F=+ \infty$.

Using Gnedenko equation \ref{gneg} we can calculate $G(t)$ as follows:

\begin{equation}
G(t) = \int_t^\infty \frac{1-F(u)}{1-F(t)} \text{d}u = \int_t^\infty \frac{ e^{-ud}}{e^{-td}}\text{d}u = \frac{1}{d} \int_{td}^\infty \frac{e^{-v}}{e^{-td}}\text{d}v= \frac{1}{d}
\label{g1_e3}
\end{equation}

According to the \citet{leadbetter} proof of Gnedenko theorem we can study both  $a_m$ and $b_m$ or $\gamma_m$  convergence as:

$$\lim_{m\to \infty} m(1-F\{ \gamma_m + xG(\gamma_m) \})=e^{-x}$$

\begin{equation}
\lim_{m \to \infty} m\Omega_d e^{-d(\gamma_m +xG(\gamma_m))} = e^{-x}
\label{g1o}
\end{equation}

then we can use the connection between $\gamma_m$ and normalising sequences to find $a_m$ and $b_m$.

By equation \ref{gamma} or using relation \ref{g1o}:

$$ \gamma_m \simeq \frac{\ln( m \Omega_d)}{d} $$

so that:

$$ a_m= d \qquad b_m=\frac{1}{d}\ln(m) + \frac{\ln(\Omega_d)}{d}$$

\paragraph{Case 2: $\mathbf{g_2}$(x)=dist(x,$\mathbf{\zeta)^{-1/\alpha}}$.} 

We can proceed as for $g_1$: 

\begin{equation}
\begin{split}
1-F(u)&= 1 - \nu(\mbox{dist}(x, \zeta)^{-1/\alpha} \leq u) \\
      &= 1 - \nu( \mbox{dist}(x, \zeta) \geq u^{-\alpha})\\
      &=  \nu( B_{u^{-\alpha}}(\zeta)) = \Omega_d u^{-\alpha d}\\
\end{split}
\end{equation}

in this case $u_F=+ \infty$.\\

\begin{equation}
\gamma_m=F^{-1}(1 -1/m)= (m\Omega_d)^{1/(\alpha d)}
\end{equation}

and, as discussed in section \ref{background}, using \citet{beirlant} choice of normalising sequences we expect: 

$$b_m=c\cdot m^{-\xi}$$

where $c \in \mathbb{R}$ is a constant.

\paragraph{Case 3: $\mathbf{g_3}$(x)=C-dist(x,$\mathbf{\zeta)^{1/\alpha}}$.}

Eventually we compute $a_m$ and $b_m$ for the $g_3$ observable class:

\begin{equation}
\begin{split}
1-F(u)&= 1 - \nu(C- \mbox{dist}(x, \zeta)^{1/\alpha} \leq u) \\
      &= 1 - \nu( \mbox{dist}(x, \zeta) \geq (C - u)^{\alpha})\\
      &=  \nu( B_{(C - u)^{\alpha}}(\zeta)) = \Omega_d(C- u)^{\alpha d}\\
\end{split}
\end{equation}

in this case $u_F= C$.\\

\begin{equation}
\gamma_m=F^{-1}(1 -1/m)= C -(m\Omega_d)^{-1/(\alpha d)}
\end{equation}

For type 3 distribution:

\begin{equation}
a_m=(u_F-\gamma_m)^{-1}, \quad b_m=u_F;
\end{equation}

\subsection{Numerical Experiments}
Since we want to show that unnormalised data may be  fitted by using the GEV  distribution $F_G(x; \mu,\sigma, \xi)$
we expect to find the following equivalence:

$$ a_m= 1/\sigma \qquad b_m=\mu$$

where, clearly, $\mu=\mu(m)$ and $\sigma=\sigma(m)$. This fact can be seen as a linear change of variable: the variable $y=a_m( x - b_m)$ has a GEV distribution  $F_G(y; \mu=0,\sigma=1, \xi)$ (that is an EV one parameter distribution with $a_m$ and $b_m$ normalising sequences) while $x$ is GEV distributed  $F_G(x; \mu=b_m,\sigma=1/a_m, \xi)$.

As we said above we now apply the previous considerations to two maps which enjoy extreme values laws and have constant density: we summarize below  the theoretical results we obtained for all three type of observables. We have obtained the results in terms of $m$ but, since we fix $k=n\cdot m$, the previous results can be translated in terms of $n$ as follows:

For $g_1$ type observable:

\begin{equation}
\sigma= \frac{1}{d} \qquad \mu \propto \frac{1}{d}\ln(k/n) 
\label{g1res}
\end{equation}

For $g_2$ type observable:

\begin{equation}
\sigma\propto n^{-1/(\alpha d)} \qquad \mu \propto n^{-1/(\alpha d)}
\label{g2res}
\end{equation}

For $g_3$ type observable:

\begin{equation}
\sigma\propto n^{1/(\alpha d)} \qquad \mu = C 
\label{g3res}
\end{equation}

Following \citet{freitas} we obtain the expression for the shape parameters: $\xi=0$ for $g_1$ type ,  $ \xi=1/(\alpha d)$ for $g_2$ type and $ \xi=-1/(\alpha d)$ for $g_3$ type.

In order to provide a numerical test of our results we consider a one-dimensional and a two dimensional map. The one dimensional
map used is a Bernoulli Shift map:

\begin{equation}
x_{t+1} = q x_t \mod 1 \qquad q>1 \in \mathbb{N}
\label{3xmod1eq}
\end{equation}

with $q=3$.\\

The considered two dimensional map is the famous  Arnold's cat map defined on the 2-torus by:

\begin{equation}
\begin{bmatrix} x_{t+1} \\ y_{t+1} \end{bmatrix} = \begin{bmatrix} 2 & 1 \\ 1 & 1 \end{bmatrix} \begin{bmatrix} x_t \\ y_t \end{bmatrix} \mod 1
\label{cat}
\end{equation}

A wide description of properties of these maps can be found in \citet{arnoldavez} and \citet{hasselblatt}.

\begin{figure}[H] \begin{center}
\advance\leftskip-1.5cm
		\includegraphics[width=1.2\textwidth]{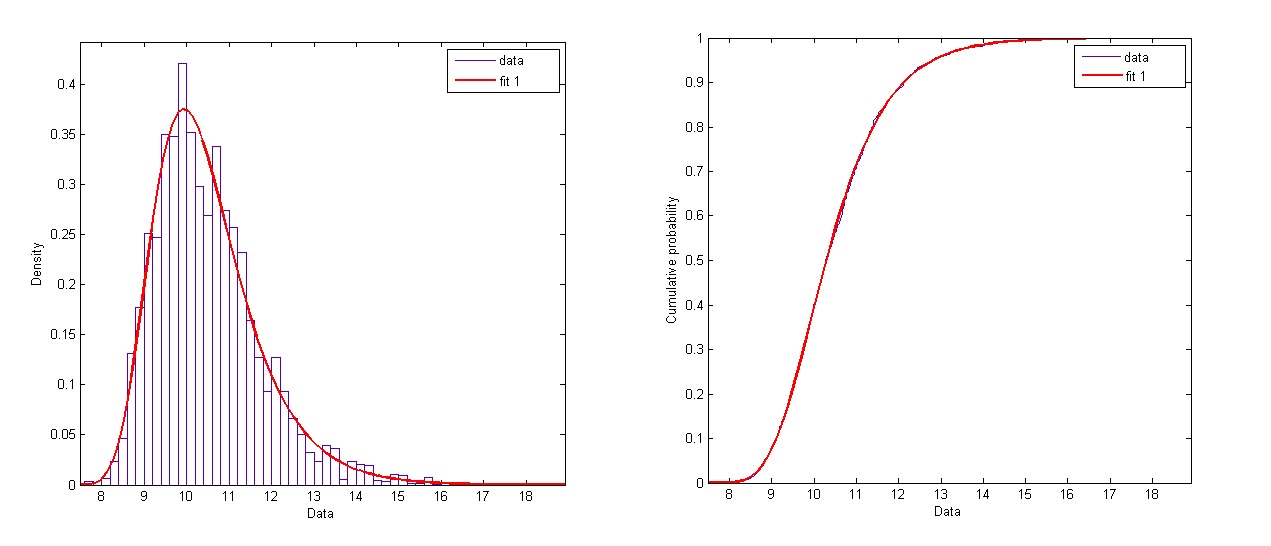}
		\caption{  Left: $g_1$ observable empirical  histogram and fitted GEV pdf. Right: $g_1$ observable empirical cdf and fitted GEV cdf. Logistic map, $n=10^4$, $m=10^4$}
	  \label{histi}
	\end{center} \end{figure}

We proceed as follows. For each map we run a long simulation up to $k$ iterations starting from a given initial condition $\zeta$. Note that 
the results - as we tested - do not depend on the choice of $\zeta$. From the trajectory we compute the sequence of observables $g_1$, $g_2$, $g_3$ as
follows dividing it into $n$ bins  each containing $m=k/n$ observations. Then, we test the degree of agreement between
the empirical distribution of the maxima and the GEV distribution according to the theoretical values presented above. A priori, it is reasonable to assume GEV as a suitable family of statistical models. For some selected values of $n$, the maxima are normalised and fitted to GEV distributions $F_G(x; \mu, \sigma, \xi)$ using a maximum likelihood method which selects values of the model parameters that produce the distribution most likely to have resulted in the observed data.\\
All the numerical analysis contained in this work has been performed using  MATLAB  Statistics Toolbox functions  such as \textit{gevfit} and \textit{gevcdf}. These functions return maximum likelihood estimates of the parameters for the generalized extreme value (GEV) distribution giving 95\% confidence intervals for estimates \citep{matlab}.\\
As in every fitting procedure, it is necessary to test the a posteriori goodness of fit. We anticipate that in every case considered, fitted distributions  passed, with maximum confidence interval, the Kolmogorov-Smirnov test described in \citet{lilliefors}. 
For illustration purposes, we present in figure \ref{histi} an empirical pdf and cdf with the corresponding fits.\\

Once $k$ is set to a given value (in our case $k=10^7$), the numerical simulations allow us to explore two limiting cases  of great interest in applications where the statistical inference is intrinsically problematic: 

\begin{enumerate}
\item $n$ is small ($m$ is large), so that we extract only few maxima, each corresponding to a \textit{very} extreme event.
\item $m$ is small ($n$ is large), so that we extract many maxima but most of those will  not be as extreme as in case 1). 
\end{enumerate}

In case 1), we have only few data - of high quality - to fit our statistical models whereas in case 2) we have many data but the sampling
may be spoiled by  the inclusion of data not giving a good representation of extreme events. We have in general that in order to obtain a reliable fit for a distribution with $p$ parameters we need $10^p$ independent data \citep{felici1} so that we expect that fit procedure gives reliable results for $n>10^3$.  As the value of $m$ determines to which extend the extracted bin maximum is representative of an extreme , below a certain value $m_{min}$
our selection procedure will be unavoidably misleading. We have no obvious theoretical argument to define the value of $m_{min}$.   We expect to obtain good fits throughout the parametric region where the constraints on $n,m$ are satisfied. Therefore, our flexibility in choosing satisfying pairs $(n,m)$ increases with larger values of $k$.\\

\begin{figure}[H] \begin{center}
\advance\leftskip-1.5cm
		\includegraphics[width=1.2\textwidth]{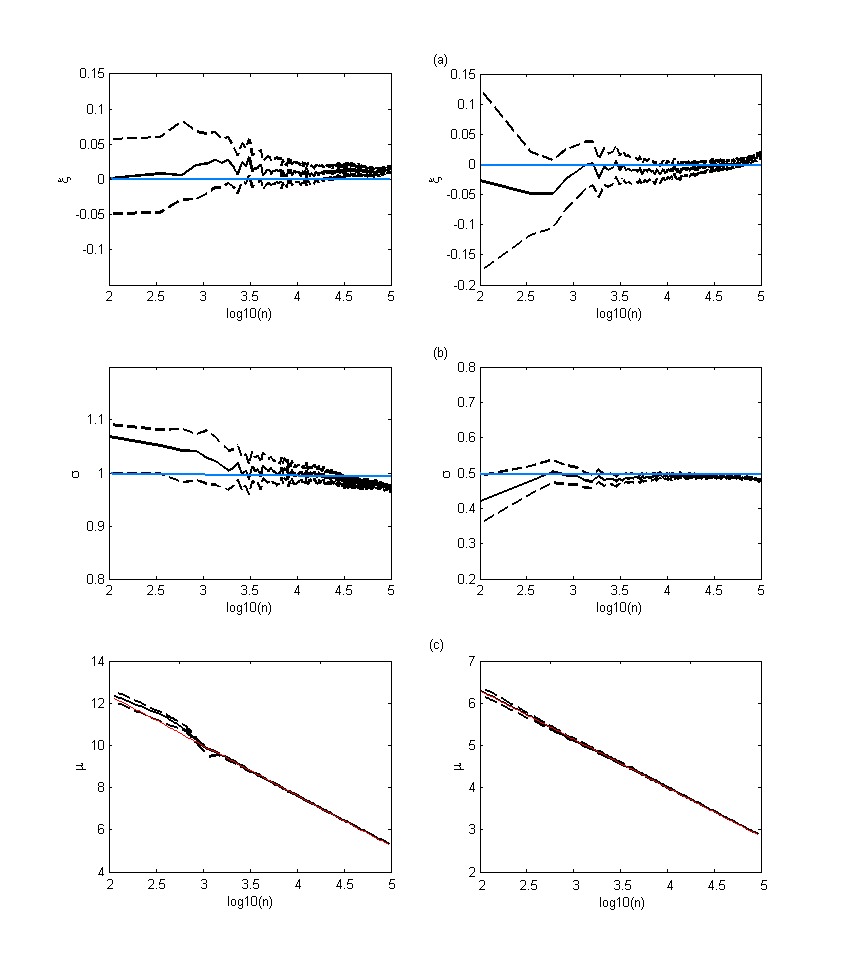}
		\caption{ $g_1$ observable, $\zeta \simeq 0.51$. \textbf{a)} $\xi$ VS $\log_{10}(n)$; \textbf{b)} $\sigma$ VS $\log_{10}(n)$; \textbf{c)} $\mu$ VS $\log(n)$. Right: Bernoulli Shift map. Left: Arnold Cat Map. Dotted lines represent computed confidence interval, gray lines represent linear fits and theoretical values.}
	  \label{g1f}
	\end{center} \end{figure}

\begin{figure}[H] \begin{center}
\advance\leftskip-1.5cm
		\includegraphics[width=1.2\textwidth]{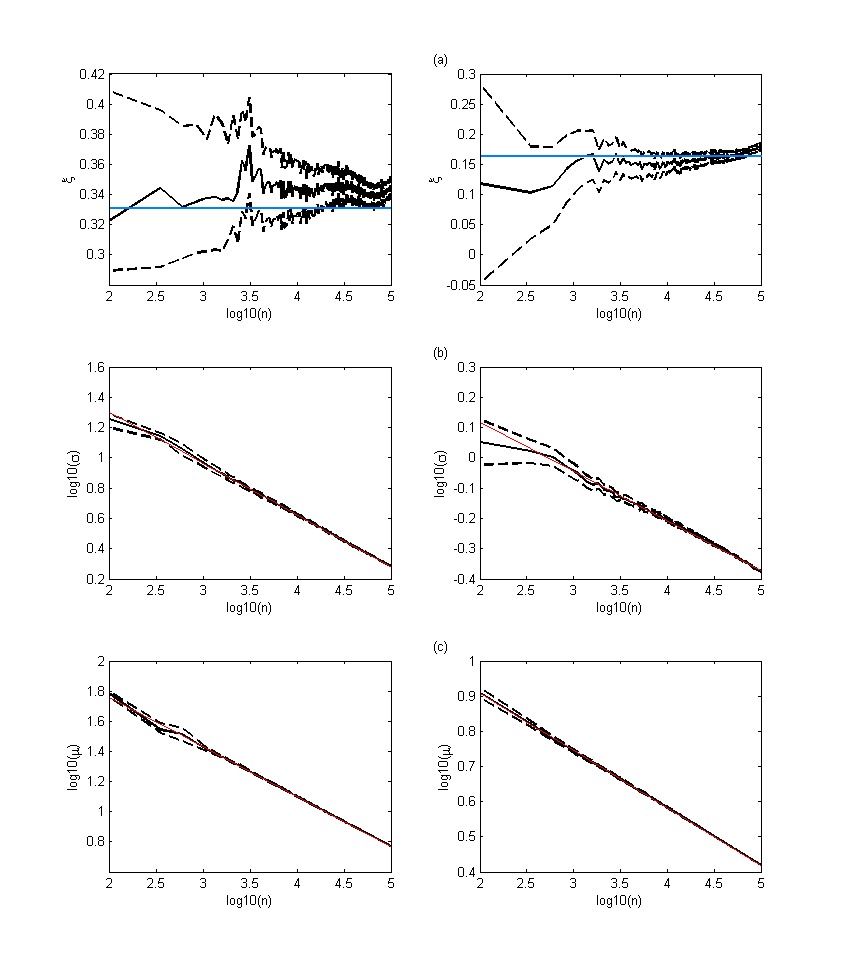}
		\caption{  $g_2$ observable, $\zeta \simeq 0.51$. \textbf{a)} $\xi$ VS $\log_{10}(n)$; \textbf{b)} $\log_{10}(\sigma)$ VS $\log_{10}(n)$; \textbf{c)} $\log_{10}(\mu)$ VS $\log_{10}(n)$. Right: Bernoulli Shift map. Left: Arnold Cat Map. Dotted lines represent computed confidence interval, gray lines represent linear fits and theoretical values.} 
	  \label{g2f}
	\end{center} \end{figure}

\begin{figure}[H] \begin{center}
\advance\leftskip-1.5cm
		\includegraphics[width=1.2\textwidth]{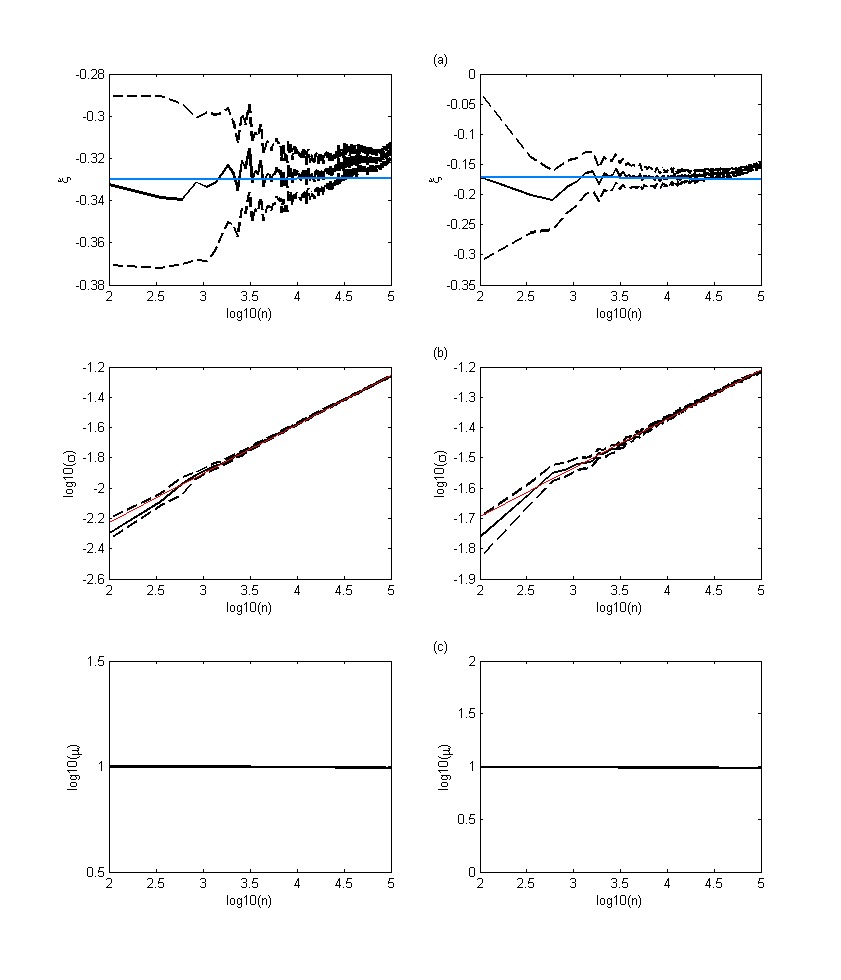}
		\caption{  $g_3$ observable, $\zeta \simeq 0.51$. \textbf{a)} $\xi$ VS $\log_{10}(n)$; \textbf{b)} $\log_{10}(\sigma)$ VS $\log_{10}(n)$; \textbf{c)} $\log_{10}(\mu)$ VS $\log_{10}(n)$. Right: Bernoulli Shift map. Left: Arnold Cat Map. Dotted lines represent computed confidence interval, gray lines represent linear fits and theoretical values.} 
	  \label{g3f}
\end{center} \end{figure}

For a $g_1$ type observable function  the behavior against $n$ of the three  parameters is presented in figure \ref{g1f}. According to equation \ref{g1res} we expect to find $\xi=0$. For relatively small values of $n$ the  sample is too small to ensure a good convergence to analytical $\xi$ and confidence intervals are wide. On the other hand we see  deviations from expected value  as  $m<10^3$ that is when $n>10^4$. For the scale parameter a similar behavior is achieved and deviations from expected theoretical values  $\sigma=1/2$ for Arnold Cat Map and $\sigma=1$ for Bernoulli Shift are found  when $n<10^3$ or $m<10^3$. Location parameter $\mu$ shows a logarithm decay with $n$ as expected from equation \ref{g1res}. A linear fit of $\mu$ in respect to $\log(n)$ is shown with agray line in figure \ref{g1f}. The linear fit computed angular coefficients $K^*$ of equation \ref{g1res} well approximate $1/d$: for Bernoulli Shift map we obtain $|K^*| = 1.001 \pm 0.001$ while for Arnold Cat map $|K^*|= 0.489 \pm 0.001$. 
We find  that $\xi$ values have  best matching with  theoretical ones with reliable confidence interval when both $n>10^3$ and $m>10^3$. These results are confirmed even for $g_2$ type and $g_3$ type observable functions as shown in figures \ref{g2f}a) and \ref{g3f}a) respectively.
We present the fit results for $\alpha=3$ but we have done tests for different $\alpha$ and for fixed $n$ and different $\alpha$.\\
For $g_2$ observable function we can also check that $\mu$ and $\sigma$ parameters follow a power law as described in eq. \ref{g2res}. In the log-log plot in Figure \ref{g2f}b), \ref{g2f}c), we can see a very clear linear
behavior. For the Bernoulli Shift map we obtain $|K^*|= 0.330 \pm 0.001 $ for $\mu$ series ,  $|K^*|= 0.341 \pm 0.001 $    for $\sigma$
in  good agreement with theoretical value of $1/3$. For Arnold Cat map we expect to find $K^*=1/6$, from the experimental data we obtain  $|K^*|= 0.163  \pm 0.001 $ for $\mu$   and  $|K^*|= 0.164 \pm 0.001 $    for $\sigma$.\\
Eventually, computing $g_3$ as observable function we expect to find  a constant value for $\mu$  while $\sigma$ has to grow with a power law in respect to $n$ as expected from
equation \ref{g3res}. 
As in $g_2$ case we expect $|K^*|=1/(\alpha d)$ and numerical results shown in figure \ref{g3f}b), \ref{g3f}c) are consistent with the theoretical one since $|K^*|=0.323 \pm 0.006$ for Bernoulli shift map and  $|K^*|= 0.162 \pm 0.006$ for Arnold Cat map.\\
 
In all cases considered the analytical behavior described in equation \ref{g2res} and \ref{g3res} is achieved and the fit quality improves if $n>10^3$ and $m>10^3$. The $g_3$ type observable constant has been chosen $C=10$. The nature of these lower bound is quite different: 

\section{Distributions of Extremes in mixing map with non-constant density measure}

\subsection{Asymptotic sequences}

The main problem when dealing with maps that have absolutely continuous but non-constant density measure $\rho(\zeta)$ is in the computation
of the integral:

\begin{equation}
\nu( B_{\delta}(\zeta)) = \int_{B_{\delta}(\zeta)} \rho(x) dx
\label{merad}
\end{equation}
where $B_{\delta}(\zeta)$ is the $d$-dimensional ball of radius $\delta$ centered in $\zeta$.\\
We have to know the value of this integral in order to evaluate $F(u)$ and, therefore, the sequences $a_m$ and $b_m$.\\
As shown in the previous section $\delta$ is linked to the observable type: in all cases, since we substitute $u=1 -1/m$,  $\delta \to 0$ means that we are interested in $m \to \infty$. \\
In this limit, a first order approximation of the previous integral is:

\begin{equation}
\nu( B_{\delta}(\zeta)) \simeq \rho(\zeta) \delta^d + \mathcal{O}(\delta^{d+1})
\label{merad2}
\end{equation}

that is valid if we are not in a neighborhood of a singular point of $\rho(\zeta)$.\\

As an example we compute the asymptotic sequences for a logistic map:

\begin{equation}
x_{t+1}=rx_t(1-x_t)
\label{logistic}
\end{equation}

with $r=4$.  This map satisfies hypothesis described in the analysis performed for Benedicks-Carleson maps in \citet{moreira}. \\

 For this map the density of the absolutely continuous invariant measure is explicit and  reads:
\begin{equation}
\rho(\zeta)=\frac{1}{\pi \sqrt{\zeta(1-\zeta)}} \qquad \zeta \in (0,1)
\label{logisticden}
\end{equation}

So that: 

 \begin{equation}
 \int_{B_{\delta}(\zeta)} \rho(\zeta) d\zeta = \frac{2}{\pi}\left[ \arcsin(\sqrt{\zeta+\delta} -\arcsin(\sqrt{\zeta-\delta} \right]
 \label{lim}
 \end{equation}

where $\zeta+\delta <1$ and $ \zeta -\delta > 0$. Since Extreme  Value Theory effectively works only if $n,m$ are large enough, the results in eq. \ref{lim} can be
replaced by  a series expansion for $\delta \to 0$:

\begin{equation}
 \frac{2}{\pi}\left[ \arcsin(\sqrt{\zeta+\delta} -\arcsin(\sqrt{\zeta-\delta} \right]= \frac{1}{\pi} \frac{2\delta}{\sqrt{\zeta(1-\zeta)}} \left[ 1 + \delta^2 P(\zeta) + ... \right]
 \label{lim2}
 \end{equation}
up to order $\delta^3$, where:

\begin{equation}
P(\zeta)= \frac{1}{8\zeta^2} - \frac{2}{\zeta(1-\zeta)} + \frac{2}{\zeta^2(1-\zeta)} + \frac{6}{\zeta^2(1-\zeta)^2}
 \label{Pa}
 \end{equation}

Using the last two equations we are able to compute asymptotic normalising sequences $a_m$ and $b_m$ for all $g_i$ observables.

\paragraph{ Case 1: $\mathbf{g_1}$(x)= -log(dist(x,$\mathbf{\zeta}$)).} 

For $g_1$ observable  functions we set $\delta=e^{-ud}$. In case of logistic map $d=1$. First we have to compute $G(t)$ using equation \ref{g1_e3} and  the expansion in eq. \label{lim3}:

\begin{equation}
G(t)= \frac{\int_t^\infty du(e^{-u} + e^{-3u}P(\zeta)}{e^{-t} + e^{-3t}P(\zeta)} \simeq 1 - \frac{2}{3} e^{-2t} P(\zeta)
\label{gt}
\end{equation}

We can compute $\gamma_m$, if $m>>1$, as follows:

\begin{equation}
 F(\gamma_m) \simeq 1 - \frac{1}{m}
\label{g1_gamman1}
\end{equation}

At the first order in eq. \ref{lim2} we get

\begin{equation}
\frac{1}{m} \simeq \frac{1}{\pi} \frac{2 e^{-\gamma_m}}{\sqrt{\zeta(1-\zeta)}}
\label{g1_gamman2}
\end{equation}

so that:

\begin{equation}
\gamma_m \simeq \ln(m)  + \ln \left(\frac{2}{\pi\sqrt{\zeta(1-\zeta)}} \right)
\label{g1_gamman3}
\end{equation}

Therefore, the sequences $a_m$ and $b_m$ if $m>>1$ are:

\begin{equation}
a_m\simeq [G(\gamma_m)]^{-1} \simeq 1 + \frac{2}{3}\frac{\pi^2}{4m^2}\zeta(\zeta-1)P(\zeta)
\label{seq1bis}
\end{equation}

\begin{equation}
b_m \simeq \gamma_m \simeq \ln(m)  + \ln \left( 2 \rho(\zeta) \right)
\label{seq2bis}
\end{equation}

\paragraph{Case 2: $\mathbf{g_2}$(x)=dist(x,$\mathbf{\zeta)^{-1/\alpha}}$.} 

We can proceed as for $g_1$ setting $\delta=(\alpha u)^{-\alpha}$, computing $\gamma_m$ we get at the first order in eq. \ref{lim2}:

\begin{equation}
\frac{1}{m} \simeq \frac{1}{\pi} \frac{2 \gamma_m^{-\alpha}}{\sqrt{\zeta(1-\zeta)}} = 2 \rho(\zeta) (\alpha\gamma_m)^{-\alpha}
\label{g2_gamman1}
\end{equation}

\begin{equation}
\gamma_m =\frac{1}{\alpha} \left(\frac{1}{2m\rho(\zeta)} \right)^{-1/\alpha}
\label{g2_gamman2}
\end{equation}

We can respectively compute $a_m$ and $b_m$ as:

\begin{equation}
a_m=\gamma_m^{-1} \qquad b_m=(2m \rho(\zeta))^{-\xi}
\label{seq2bis2}
\end{equation}

\paragraph{Case 3: $\mathbf{g_3}$(x)=C-dist(x,$\mathbf{\zeta)^{1/\alpha}}$.}  

As in the previous cases, we compute $\gamma_m$ up to the first order setting $\delta=[\alpha(C-\gamma_m)]^\alpha$:

\begin{equation}
\frac{1}{m} \simeq \frac{1}{\pi} \frac{2 [\alpha(C-\gamma_m)]^{\alpha}}{\sqrt{\zeta(1-\zeta)}} = 2 \rho(\zeta)[\alpha(C-\gamma_m)]^{\alpha}
\label{g3_gamman1}
\end{equation}

\begin{equation}
\gamma_m = C - \frac{1}{\alpha} \left(\frac{1}{2m\rho(\zeta)} \right)^{1/\alpha}
\label{g3_gamman2}
\end{equation}

For type 3 distribution:

\begin{equation}
a_m=(u_F-\gamma_m)^{-1}, \quad b_m=u_F;
\label{seq3bis}
\end{equation}

where $u_f=C$.

\subsection{Numerical experiment on the logistic map}

Following the same procedure detailed in section 3.2, we want to show the equivalence between EV computed normalising sequences $a_m$ and $b_m$ and the parameters of a GEV distribution obtained directly fitting the data even in case of logistic map that has not constant density measure. Using eq. \ref{seq1bis}-\ref{seq2bis} for $g_1$, we obtain the following theoretical expression:

\begin{equation}
\sigma(m,\zeta) \simeq 1 + \frac{2}{3}\frac{\pi^2}{4m^2}\zeta(\zeta-1)P(\zeta) \qquad \mu(m,\zeta) \simeq \ln(m) + \ln(2\rho(\zeta))
\label{g1log}
\end{equation}

From eq. \ref{seq2bis}, for  $g_2$ observable type, we write:

\begin{equation}
\sigma(m,\zeta) \simeq \frac{1}{\alpha} (2m\rho(\zeta))^{1\over {\alpha}} \qquad   \mu(m,\zeta) \simeq  (2m\rho(\zeta))^{1\over  {\alpha}}
\label{g2log}
\end{equation}

and in $g_3$ case  using eq. \ref{seq3bis}, we expect to find:
\begin{equation}
\sigma(m,\zeta) \simeq \frac{1}{\alpha}(2m\rho(\zeta))^{-1/\alpha} \qquad   \mu(m,\zeta) \simeq C = u_F
\label{g3log}
\end{equation}

Values of $\xi$ are independent on density and, as stated in Freitas' $\xi=0$ for $g_1$ type ,  $ \xi=1/(\alpha d)$ for $g_2$ type and $ \xi=-1/(\alpha d)$ for $g_3$ type.\\

\begin{figure}[H] \begin{center}
\advance\leftskip-1.5cm
		\includegraphics[width=1.0\textwidth]{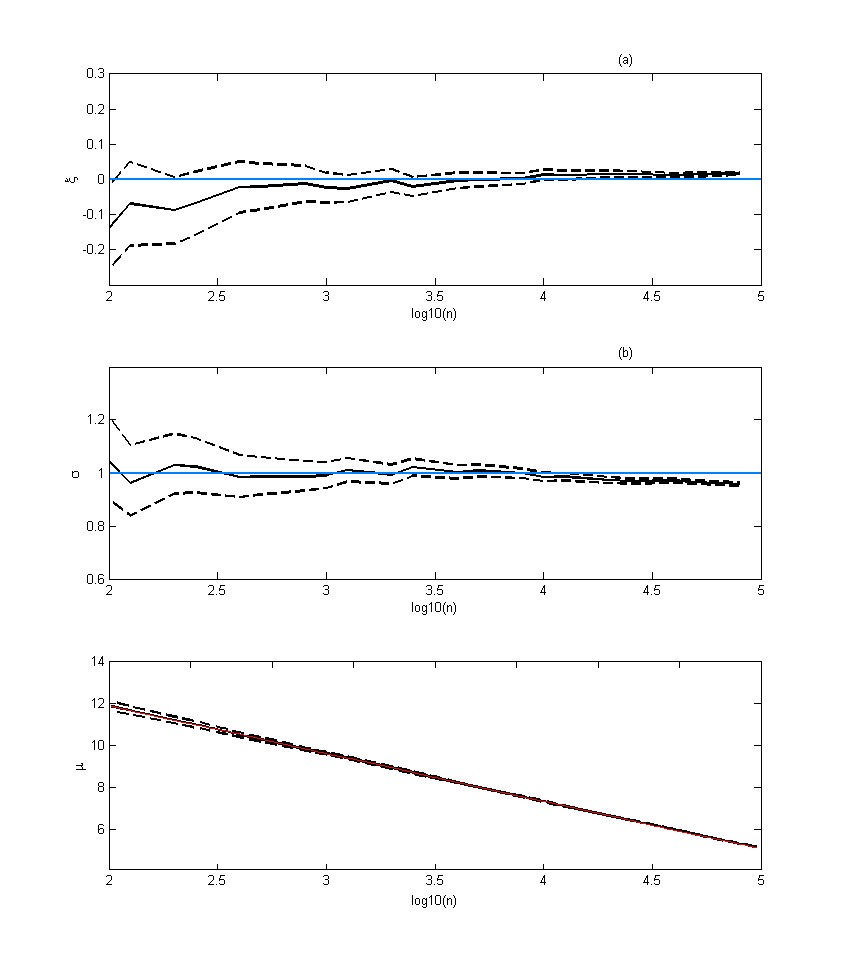}
		\caption{ $g_1$ observable, $\zeta = 0.31$. \textbf{a)} $\xi$ VS $\log_{10}(n)$; \textbf{b)} $\sigma$ VS $\log_{10}(n)$; \textbf{c)} $\mu$ VS $\log(n)$. Logistic map. Dotted lines represent computed confidence interval, gray lines represent linear fits and theoretical values.}
	  \label{l1}
	\end{center} \end{figure}

\begin{figure}[H] \begin{center}
\advance\leftskip-1.5cm
		\includegraphics[width=1.2\textwidth]{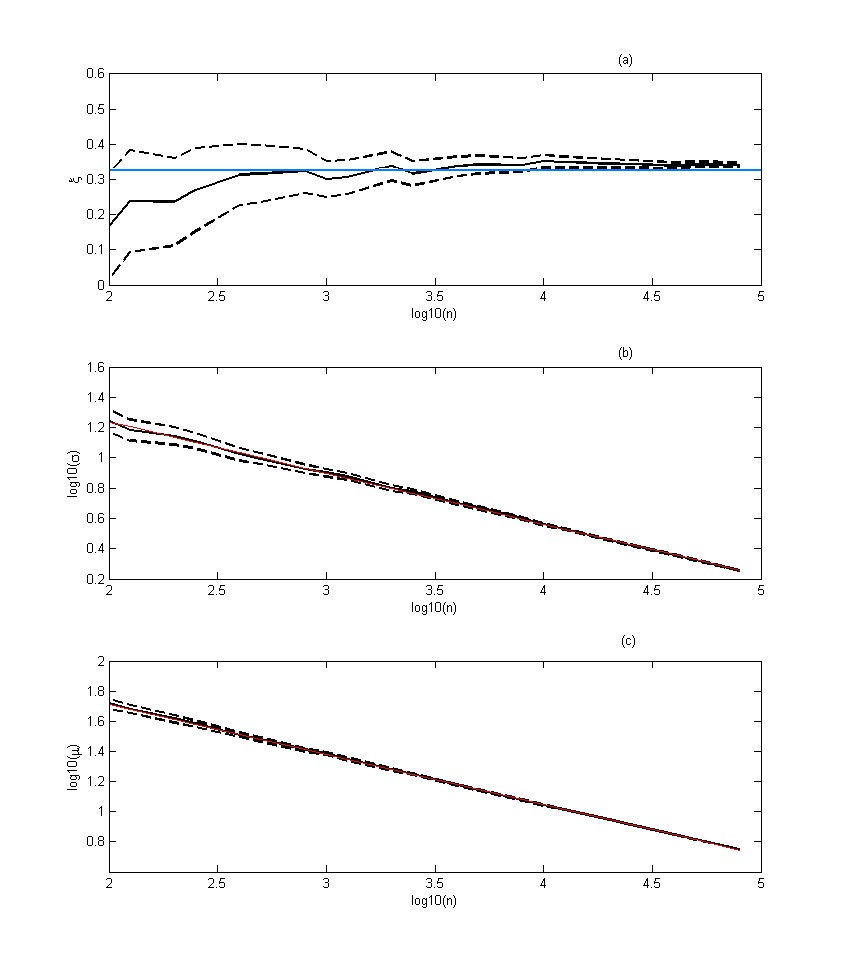}
		\caption{  $g_2$ observable, $\zeta = 0.3$. \textbf{a)} $\xi$ VS $\log_{10}(n)$; \textbf{b)} $\log_{10}(\sigma)$ VS $\log_{10}(n)$; \textbf{c)} $\log_{10}(\mu)$ VS $\log_{10}(n)$. Logistic map. Dotted lines represent computed confidence interval, gray lines represent linear fits and theoretical values.} 
	  \label{l2}
	\end{center} \end{figure}

\begin{figure}[H] \begin{center}
\advance\leftskip-1.5cm
		\includegraphics[width=1.2\textwidth]{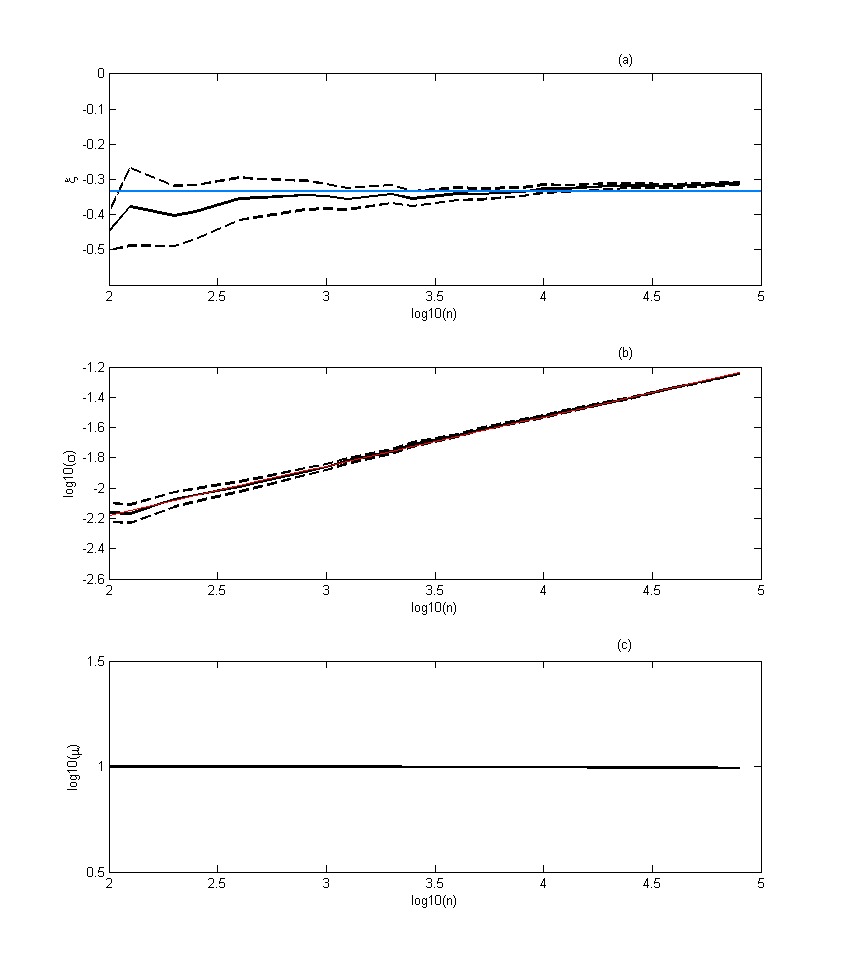}
		\caption{  $g_3$ observable, $\zeta= 0.3$. \textbf{a)} $\xi$ VS $\log_{10}(n)$; \textbf{b)} $\log_{10}(\sigma)$ VS $\log_{10}(n)$; \textbf{c)} $\log_{10}(\mu)$ VS $\log_{10}(n)$. Logistic map. Dotted lines represent computed confidence interval, gray lines represent  linear fits and theoretical values.} 
	  \label{l3}
	\end{center} \end{figure}

In figures \ref{l1}-\ref{l3} we presents a numerical test of the asymptotic behavior described in equations \ref{g1log} - \ref{g3log} on logistic map for $d=1$ , $a=3$, $C=u_F=10$, $\zeta=0.3$ against the variable $n$.  As shown in previous section,  block maxima approach works well with maps with constant density measure when $ n $ and $m$ are at least $10^3$: In fact, regarding $\xi$ parameter. Significant deviations from the theoretical value are achieved when $n<1000$ or $m<1000$ even in the case of the Logistic Map.\\
Regarding $\mu$ and $\sigma$, for $g_1$ observable a linear fit of $\mu$ in respect to $\log(n)$ give us $|K^*|=0.999\pm0.002$, while $\sigma$ shows the same behavior of $\xi$ since the best agreement with theoretical value $\sigma=1$ is achieved when $n,m>10^3$. In the log-log plots of figure \ref{l2}b), \ref{l2}c) for $g_2$ observable, we can observe again the expected linear
behavior for $\mu$ and $\sigma$ with $|K^*|$ corresponding to $1/(\alpha d)$. From numerical fit we obtain $|K^*|= 0.3334 \pm 0.0007 $ for $\mu$ series and  $|K^*|= 0.337 \pm 0.002 $    for $\sigma$ in  good agreement with theoretical value of $1/3$. By applying a linear fit to the log-log plot in figure \ref{l3}b), the angular coefficient corresponding to $\sigma$ series is $|K|= 0.323 \pm 0.003$ again consistent with the theory.\\

For a logistic map we can also check the GEV behavior in respect to initial conditions.  If we fix $n^*=m^*=10^3$ and fit our data to GEV distribution for $10^3$ different $\zeta \in (0,1)$ an asymptotic behavior is reached as shown from the previous analysis.
For $g_1$ observable function  we have observed  that the first order approximation works well for all three parameters. Deviation from this behavior are achieved for $\zeta \to 1$ and $\zeta \to 0$ as the measure become singular when we move to these points and we should take in account other terms of the series expansion. Numerically, we found that deviations from first order approximation are meaningful only if $\zeta <10^{-3}$ and $\zeta> 1-10^{-3}$.  Averaging over $\zeta$ both $\xi$ and $\sigma$ we obtain $<\xi> = 1.000 \pm 0.009 $ and  $<\sigma>=1.00 \pm 0.03$ where the uncertainties are computed with respect to the estimator. Since we expect $\xi=0$ and $\sigma=1$ at zero order approximation,  numerical results are consistent with the theoretical ones; furthermore, experimental data are normally distributed around theoretical values.\\
Asymptotic expansion also works well for $g_2$ observables: we obtain $<\xi> = 0.334 \pm 0.001 $ in excellent agreement with theoretical value $\xi=1/3$.
Eventually, in $g_3$,  averaging $\xi$ over different initial conditions we get  $<\xi> = -0.334 \pm 0.002$ that is again consistent to theoretical value -1/3.

%$\sigma(\zeta)$ and $\xi(\zeta)$ values are fluctuating around the theoretical values (1 and 0 respectively): averaging $\xi$ over $\zeta$ we obtain $  (1.0 \pm 0.9) \cdot 10^{-3} $, while doing the same with $\sigma$ gives a value of $1.00 \pm 0.03$ .\\
%
%The case of $g_2$ observable  function is shown in figure \ref{l2}:  both $\mu(\zeta)/[2n^*\rho(\zeta)]^{1/3}$ and $\sigma(\zeta)/[(1/3)(2n^*\rho(\zeta))^{1/3}]$ show a constant behavior around 1 as expected and averaging for all initial conditions $\xi$ we obtain $0.334 \pm 0.001$ which is consistent with the theoretical value 1/3.\\
%
%Eventually, in  $g_3$ case shown in figure \ref{l3}, $\mu(\zeta)$ is almost constant around 10 that is its expected theoretical value from eq. \ref{g3log} and $\sigma(\zeta)/[(1/3)(2n^*\rho(\zeta))^{-1/3}]$ against $\zeta$ is almost constant around 1. Averaging $\xi$ over different initial conditions we get  $-0.334 \pm 0.002$ that is again consistent to theoretical value -1/3.

\section{The case of regular maps}
\citet{freitas2008} have posed the problem of dependent extreme values in dynamical systems that show uniform quasi periodic motion. Here we try to investigate this problem numerically. We  have  used  a one-dimensional and a bi-dimensional discrete map. The first one is the irrational translation on the torus defined by:

\begin{equation}
x_{t+1} =  x_t + \beta \mod 1 \qquad \beta \in [0,1]  \setminus \mathbb{Q}
\label{torus_trasl}
\end{equation}

And for the bidimensional case, we use the so called standard map:

\begin{equation} 
  y_{t+1}=y_t+\frac{\lambda}{2\pi} \sin(2\pi x_t) \mod 1; 
  \qquad x_{t+1}=x_t+y_{t+1}\mod 1. 
  \label{stdmap}
\end{equation}

with $\lambda=10^{-4}$. For this value of $\lambda$, the standard map exhibits a regular behavior and it is not mixing, as well as torus translations. This means that these maps fail in satisfying hypothesis $D_2$ and $D'$ and moreover they do not enjoy as well an exponential hitting time statistics. About this latter statistics, it is however known that  it exists for torus translation and it is given by a  particular piecewise linear function or a  uniform distribution depending on which sequence of sets $A_k$ is considered  \citep{coelho1996limit}. In a similar way, a non-exponential Hitting Time Statistics (HTS) is achieved for standard map when $\lambda<<1$ as well as for a skew map, that is a standard map with $\lambda=0$ \citep{buric2005statistics}. Therefore we expect not to obtain a GEV distribution of any type using $g_i$ observables.  \\

We have pointed out that the observable functions choice is  crucial in order to observe some kind of distribution of extreme values when we are dealing with dynamical systems instead of stochastic series. \citet{nicolis}  have shown how it is possible to obtain an analytical EV distribution which does not belong to GEV family choosing a simple observable: they considered the series  of  distances  between the iterated trajectory and  the initial condition. Using the same notation of section \ref{background} we can write:

$$ Y_m(x=f^t\zeta)= \mbox{dist}(f^t\zeta, \zeta) \qquad \hat{M}_m=\min \{Y_0, ... Y_{m-1} \}$$

For this observable they have shown that the cumulative distribution $F(x)=P\{a_m(\hat{M}_m-b_m) \leq x\}$  of a uniform quasi periodic motions is not smooth but piecewise linear (\citet{nicolis}, figure 3). Furthermore  slop changes of $F(x)$ can be explained by constructing the intersections between  different
iterates of equation \ref{torus_trasl}. 
$F(x)$ must correspond  to a density distribution  continuous obtained as a composition of box functions:  each box must be related to a change in the slope of $F(x)$. \\

The numerical results we report below confirm that for the maps \ref{torus_trasl} and \ref{stdmap}
the distributions of maxima for various observables cannot be fitted
with a GEV since  they are multi modal. We recall that the   return  times into a
sphere of vanishing radius do not have a spectrum,   if the orbits
have the same frequency,  whereas a spectrum  appears if the frequency
 varies continuously with the action,  as in the standard map 
for $\lambda$ close to zero \citep{turchetti2}.  Since  the EV statistics refers to a
single orbit, no change due to the local mixing,  which  insures  the
existence of a return times  spectrum \citep{turchetti2},  can be observed.
Considering that  the GEV exists when the system is mixing and does
not when it is integrable, one might use the quality of fit to GEV as
a dynamical indicator, for systems  which exhibit regions with
different dynamical properties,  ranging from integrable to  mixing as
it occurs for the standard map when $\lambda$ is order 1.
Indeed we expect that in the neighborhood of a low order resonance,
where the omoclinic tangle of intersecting separatrices appears, a 
GEV fit is possible. Preliminary computations carried out for the
standard map and for a model with parametric resonance confirm this
claim, that will be carefully tested in the near future.\\

Using the theoretical framework provided in \citet{nicolis} we check numerically the behavior of maps described in eq. \ref{torus_trasl}-\ref{stdmap} analysing EV distributions for $g_i$ observable functions.
Proceeding as in section 3 for mixing maps, we try to perform a fit to GEV distribution starting  with different initial conditions $\zeta$, a set of different $\alpha$ values and $(n,m)$ combinations. In all cases analysed  the Kolmogorov Smirnov test  fails and this means that GEV distribution is not useful to describe the behavior of this kind of statistics. This result is in agreement with \citet{freitasNuovo} but we may find out which kind of empirical distribution is obtained.\\

\begin{figure}[H] \begin{center}
\advance\leftskip-1.5cm
		\includegraphics[width=1.2\textwidth]{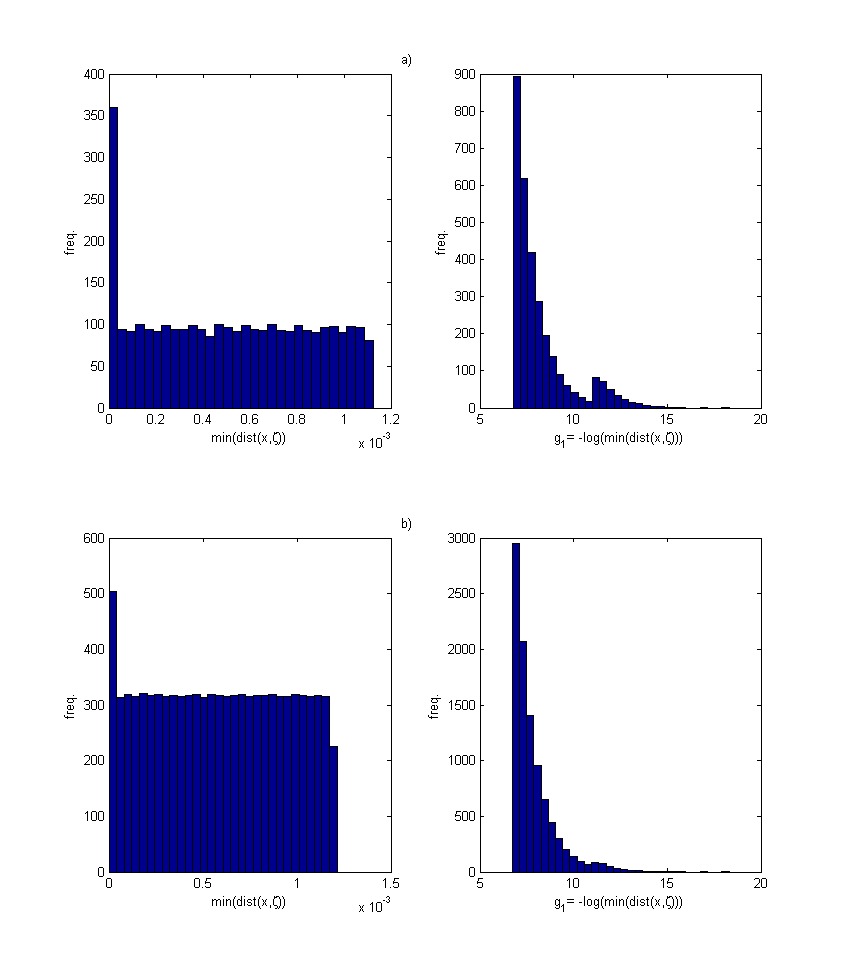}
		\caption{ Histogram of maxima for $g_1$ type observable function, standard map, $x_0=y_0= \sqrt{2}-1$.
		Left: series of $\min(\mbox{dist}(f^t\zeta, \zeta))$. Right: series of $g_1=-\log(\min(\mbox{dist}(f^t\zeta, \zeta)))$. a) $n=3300$, $m=3300$. b) $n=10000$, $m=1000$.}
	  \label{hist}
	\end{center} \end{figure}

Looking in details at $M_m$ histograms that correspond to empirical density distributions, they appear always to be multi modal and each mode have a well defined shape: for $g_1$ type observable function modes are exponential while, for $g_2$ and $g_3$, their shape depends on $\alpha$ value of observable function. Furthermore, the number of modes and their positions are highly dependent on both $n,m$ and initial conditions.\\
Using \citet{nicolis} results it is possible to understand why we obtain this kind of histograms: since density distribution  of $\hat{M}_m$ is  a composition of box functions, when we apply $g_i$ observables we  modulate it changing the shape of the boxes. Therefore, we obtain a multi modal distribution modified according to the observable functions $g_i$.\\
An example is shown in figure \ref{hist} for standard map: the left figures correspond to the histogram of the minimum distance obtained without computing $g_i$ observable and reproduce a composition of box functions. The figures in the right show how this distribution is modified by applying $g_1$ observable to the series of minimum distances.  We can see two exponential modes, while the third is hidden in the linear scale but can be highlighted using a log-scale. The upper figures are drawn using  $n=3300$, $m=3300$, the lower with $n=10000$, $m=1000$.

\section{Concluding Remarks}

EVT was developed to study a wide class of problems of great interest in different disciplines: the need of modeling events that occur with very small probability comes from the fact that they can affect in a strong way several socioeconomic activities: floods, insurance losses, earthquakes, catastrophes. A very extensive account of EVT applications has been recently given in \citet{ghil2010extreme}. EVT was applied on limited data series using the block-maxima approach  facing the problem of  having a good statistics of extreme values retaining a sufficient number of observation in each bin. Often, since no theoretical \textit{a priori} values of GEV parameters are available for this kind of applications, we may obtain a biased fit to GEV distribution even if  tests of statistical significance succeed. The recent development of an extreme value theory in dynamical systems give us the theoretical  framework to test the consistency of block-maxima approach when analytical results for distribution parameters are available. This theory relies on the global properties of the dynamical systems considered (such as the degree of mixing or the decay rate of the Hitting Time Statistics) but also on the observable functions we chose.\\

Our main finding is that a block-maxima approach for GEV distribution is totally equivalent to fit  an EV distribution after normalising sequences are computed. To prove this  we have derived analytical expressions for $a_m$ and $b_m$ normalising sequences, showing that $\mu$ and $\sigma$ of fitted GEV distribution can replace them. This approach works for maps that have an absolutely continuous invariant  measure and retain some mixing properties that can be directly related to the exponential decay of HTS. 
Since GEV approach does not require the a-priori knowledge of the measure density that is instead require by the EV approach, it is possible to use it in many numerical applications.\\

Furthermore, if we compare analytical and numerical results we can study what is the minimum number of maxima and how big the set of observations in which the maximum is taken has to be. To accomplish this goal we have analysed maps with constant density measure finding that a good agreement between numerical and analytical value is achieved when both the number of maxima $n$  and the observations per bin $m$  are at least $10^3$. We remark that  the fits have passed Kolmogorov Smirnov test with maximum confidence interval even if $n<10^3$ or $<m<10^3$ so that parametric or non parametric tests are not the only thing to take in account when dealing with extreme value distributions: if maxima are not proper extreme values (which means $m$ is not large enough) the fit is good but parameters are different from expected values. The lower bound of $n$ can be explained using the argument that a fit to a 3-parameters distribution needs at least $10^3$ independent data to give reliable informations.\\  
Therefore, we checked that in case of non-constant absolutely continuous density measure the asymptotic expressions used to compute $\mu$ and $\sigma$ works when we consider $n$ and $m$ of order $10^3$. For logistic map the numerical values of parameters  we obtain averaging over different initial conditions are totally in agreement with the theoretical ones. 
In regular maps, as expected, the fit to  a GEV distribution is unreliable. We obtain a multi modal distribution, that, for the analyzed maps, is the result of a composition of  modes in which the shape depends on observable types. This  behavior can be explained pointing out that this kind of systems have not an exponential HTS decay and therefore have no EV law for  observables of type $g_i$.\\

To conclude, we claim that we have provided a reliable way to investigate properties of
extreme values in mixing dynamical systems  which may satisfy mixing conditions (like $D_2$ and $D'$),  finding an equivalence among $a_m$, $b_m$, $\mu$ and $\sigma$ behavior for absolutely continuous measures. In our future work we intend to address the case of singular measure. Recently the theorem was generalised to the case of non smooth observations and therefore it holds also with non absolutely  continuous invariant probability measure  \citep{freitasNuovo}. In this case we expect the same for  all the procedure described here. Understanding the extreme values behavior for singular measures will be crucial to apply proficiently this analysis to operative geophysical models since in these case we are always dealing with singular measures.  In this way we will provide a complete tool to study  extreme events in complex dynamical systems   used in geophysical or financial applications.

\section{Acknowledgments}

S.V. was supported by the CNRS-PEPS Project \textit{Mathematical Methods of
Climate Models}, and he thanks the GDRE Grefi-Mefi for having supported
exchanges with Italy. V.L. and D.F. acknowledge the financial support of the
EU FP7-ERC project NAMASTE: Thermodynamics of the Climate System.

\bibliographystyle{plainnat}  
\bibliography{Biblio}

\end{document}